\documentclass[reqno,a4paper]{amsart}%{birkmult}%
\usepackage{amssymb,url,graphicx}
\usepackage[hypertex]{hyperref}%,pagebackref
\date{}
\allowdisplaybreaks[4]
\theoremstyle{plain}
\newtheorem{thm}{Theorem}
\theoremstyle{remark}
\newtheorem{rem}{Remark}
\DeclareMathOperator{\td}{d\mspace{-2mu}}

\begin{document}

\title{Sharp inequalities for the psi function and harmonic numbers}

\author[F. Qi]{Feng Qi}
\address[F. Qi]{Research Institute of Mathematical Inequality Theory, Henan Polytechnic University, Jiaozuo City, Henan Province, 454010, China}
\email{\href{mailto: F. Qi <qifeng618@gmail.com>}{qifeng618@gmail.com}, \href{mailto: F. Qi <qifeng618@hotmail.com>}{qifeng618@hotmail.com}, \href{mailto: F. Qi <qifeng618@qq.com>}{qifeng618@qq.com}}
\urladdr{\url{http://qifeng618.spaces.live.com}}

\author[B.-N. Guo]{Bai-Ni Guo}
\address[B.-N. Guo]{School of Mathematics and Informatics,
Henan Polytechnic University, Jiaozuo City, Henan Province, 454010, China} \email{\href{mailto: B.-N.
Guo <bai.ni.guo@gmail.com>}{bai.ni.guo@gmail.com}, \href{mailto: B.-N. Guo
<bai.ni.guo@hotmail.com>}{bai.ni.guo@hotmail.com}} \urladdr{\url{http://guobaini.spaces.live.com}}

\begin{abstract}
In this paper, two sharp inequalities for bounding the psi function $\psi$ and the harmonic numbers $H_n$ are established respectively, some results in [I. Muqattash and M. Yahdi, \textit{Infinite family of approximations of the Digamma function}, Math. Comput. Modelling \textbf{43} (2006), 1329\nobreakdash--1336.] are improved, and some remarks are given.
\end{abstract}

\subjclass[2000]{Primary 33B15, Secondary 26D15}

\keywords{sharp inequality, psi function, harmonic numbers, monotonicity, convex function}

\thanks{The authors were supported in part by the China Scholarship Council}

\thanks{This paper was typeset using \AmS-\LaTeX}

\maketitle

\section{Introduction}

It ie well-known that the classical Euler's gamma function is defined by
\begin{equation}
\Gamma(x)=\int_0^\infty e^{-t}t^{x-1}\td t
\end{equation}
for $x>0$ and the derivative of its logarithm is called the psi or digamma function and denoted by $\psi(x)$ for $x>0$.
\par
In \cite{Muqattash-Yahdi}, an infinite family of approximations for the psi function $\psi(x)$ on $(0,\infty)$, denoted as $\{I_a,a\in[0,1]\}$, where
\begin{equation}
  I_a(x)=\ln(x+a)-\frac1x,
\end{equation}
was constructed. Among other things, Corollary~2.3 and Theorem~3.2 in \cite{Muqattash-Yahdi} may be recited as follows:
\begin{enumerate}
  \item
  For all $x\in(0,\infty)$,
  \begin{equation}\label{corollary2.3}
    \ln(x+1)-\frac1x\ge\psi(x)\ge\ln x-\frac1x.
  \end{equation}
  \item
  For every $x\in(0,\infty)$, there exists an $a\in[0,1]$ such that $\psi(x)=I_a(x)$.
\end{enumerate}
\par
\begin{figure}[htbp]
  \includegraphics[width=0.7\textwidth]{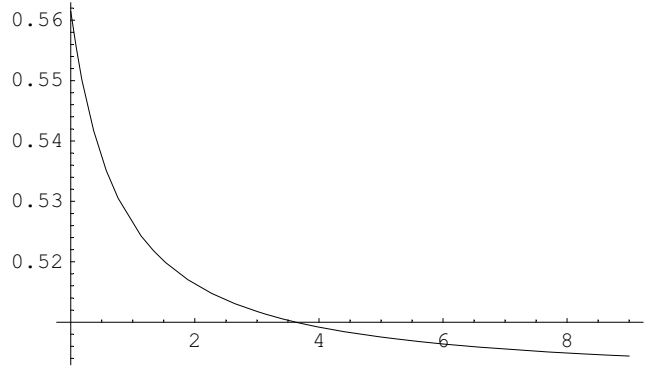}\\
  \caption{The graph on $(0,9)$ and two limits as $x\to0^+$ and $\infty$ of $Q(x)$}
  \label{Infinite-family-Digamma.tex.eps}
\end{figure}

The inequality~\eqref{corollary2.3} may be rearranged as
\begin{equation}\label{Q(x)-dfn}
1\ge\exp\biggl(\psi(x)+\frac1x\biggr)-x\triangleq Q(x)\ge0
\end{equation}
for $x\in(0,\infty)$. Since $\Gamma(x+1)=x\Gamma(x)$ for $x>0$, taking the logarithm of this recurrent formula and differentiating yields
\begin{equation}\label{psisymp4}
\psi(x+1)=\psi(x)+\frac1x.
\end{equation}
As a result, the function $Q(x)$ defined in~\eqref{Q(x)-dfn} may be rearranged as
\begin{equation}\label{Q(x)-dfn-2}
Q(x)=e^{\psi(x+1)}-x
\end{equation}
for $x\in(0,\infty)$. The graph of $Q(x)$ in the interval $(0,9)$, plotted by the famous software M\textsc{athematica}~5.2, and the limits
\begin{equation}\label{2-limits-Q(x)}
\lim_{x\to0^+}Q(x)=e^{-\gamma}\quad \text{and} \quad \lim_{x\to\infty}Q(x)=\frac12,
\end{equation}
calculated also by M\textsc{athematica}~5.2, see Figure~\ref{Infinite-family-Digamma.tex.eps},
show that the function $Q(x)$ is possibly decreasing on $(0,\infty)$ and $\frac12<Q(x)<e^{-\gamma}=0.56\dotsm$ for $x\in(0,\infty)$, where $\gamma=0.577\dotsm$ stands for the Euler-Mascheroni's constant. These analyses imply that conclusions obtained in \cite{Muqattash-Yahdi} may possibly be refined and restated more accurately.
\par
Our main results are included in the following theorems.

\begin{thm}\label{Qi-Modell-thm1}
For all $x\in(0,\infty)$,
  \begin{equation}\label{corollary2.3-rew}
    \ln\biggl(x+\frac12\biggr)-\frac1x< \psi(x)< \ln(x+e^{-\gamma})-\frac1x.
  \end{equation}
The constants $\frac12$ and $e^{-\gamma}=0.56\dotsm$ in~\eqref{corollary2.3-rew} are the best possible.
\end{thm}

\begin{thm}\label{Q(x)-CM}
The function $Q(x)$ defined in~\eqref{Q(x)-dfn-2} is strictly decreasing and strictly convex on $(-1,\infty)$.
\end{thm}

\begin{thm}\label{Qi-Modell-thm2}
For every $x\in(0,\infty)$, there exists a unique number $a\in\bigl(e^{-\gamma},\frac12\bigr)$ such that $\psi(x)=I_a(x)$. Conversely, for every $a\in\bigl(e^{-\gamma},\frac12\bigr)$, there exists a unique number $x\in(0,\infty)$ such that $\psi(x)=I_a(x)$.
\end{thm}

In \cite[Definition~3.5]{Muqattash-Yahdi}, the so-called error of the approximation $\psi(x)\approx I_a(x)$ was defined by
\begin{equation}
E_a(x)=\psi(x)-I_a(x)=\psi(x)-\ln(x+a)+\frac1x
\end{equation}
for $x>0$ and $a\in[0,1]$. In \cite[Theorem~3.7]{Muqattash-Yahdi}, it was proved that the errors $E_a(x)$ for $x\in[2,\infty)$ and $a\in[0,1]$ are uniformly bounded between $-\ln\frac32$ and $\ln\frac32$. This can be restated as follows.

\begin{thm}\label{error-estimate}
Let $c=1.4616321\dotsm$ is the only positive root of the psi function $\psi(x)$ on $(0,\infty)$. For $x\in(c,\infty)$ and $a\in\bigl(e^{-\gamma},\frac12\bigr)$, the following conclusions are valid:
\begin{enumerate}
\item
The errors $E_a(x)$ are uniformly bounded between $-\ln\frac{c+e^{-\gamma}}{c+1/2}$ and $\ln\frac{c+e^{-\gamma}}{c+1/2}$;
\item
$\psi(x)=\ln(x+a)-\frac1x+O\bigl(\ln\frac{x+e^{-\gamma}}{x+1/2}\bigr)$.
\end{enumerate}
\end{thm}

It is well-known that the $n$-th harmonic numbers are defined for $n\in\mathbb{N}$ by
\begin{equation}
H_n=\sum_{k=1}^n\frac1k
\end{equation}
and that $H_n$ can be expressed in terms of the psi function $\psi(x)$ by
\begin{equation}\label{H_n-dfn}
H_n=\psi(n+1)+\gamma.
\end{equation}
\par
By virtue of the decreasing monotonicity of the function $Q(x)$ and the formula~\eqref{H_n-dfn}, the following new bounds for $H_n$ are derived as follows.

\begin{thm}\label{harmonic-no-vu}
For $n\in\mathbb{N}$,
\begin{equation}\label{new-bounds-harmonic-No}
\ln\biggl(n+\frac12\biggr)+\gamma<H_n(n)\le\ln\bigl(n+e^{1-\gamma}-1\bigr)+\gamma.
\end{equation}
\end{thm}

In the final section, some remarks about above conclusions are given.

\section{Proofs of theorems}

\begin{proof}[Proof of Theorem~\ref{Qi-Modell-thm1}]
Since $\psi(1)=-\gamma$, the first limit in~\eqref{2-limits-Q(x)} is valid clearly.
\par
In \cite{qi-cui-jmaa}, it was derived that
\begin{equation}\label{lem1}
\frac{1}{2x}-\frac{1}{12x^2}<\psi(x+1)-\ln x<\frac{1}{2x}
\end{equation}
for $x>0$. Then
\begin{equation*}
 xe^{1/(2x)-1/(12x^2)}-x <Q(x)<xe^{1/(2x)}-x.
\end{equation*}
It is easy to check that both bounds for $Q(x)$ above tend to $\frac12$ as $x\to\infty$. The second limit in~\eqref{2-limits-Q(x)} follows.
\par
In~\cite[Lemma~1.1]{batir-new} and~\cite[Lemma~1.1]{batir-new-rgmia}, the inequality
\begin{equation}\label{batir-alzer-ineq}
\psi'(x)e^{\psi(x)}<1
\end{equation}
for $x>0$ was obtained. This means that
\begin{equation}\label{Q(x)-1-der}
  [Q(x)]'=\psi'(x+1)e^{\psi(x+1)}-1<0
\end{equation}
for $x>-1$. Consequently, the function $Q(x)$ is strictly decreasing on $(-1,\infty)$.
\par
Combining the decreasing monotonicity with two limits in~\eqref{2-limits-Q(x)} of $Q(x)$ leads obviously to Theorem~\ref{Qi-Modell-thm1}.
\end{proof}

\begin{proof}[Proof of Theorem~\ref{Q(x)-CM}]
The decreasing monotonicity has been proved in the proof of Theorem~\ref{Qi-Modell-thm1}.
\par
Differentiating~\eqref{Q(x)-1-der} once again gives
\begin{equation*}
[Q(x)]''=\bigl\{\psi''(x+1)+[\psi'(x+1)]^2\bigr\}e^{\psi(x+1)}.
\end{equation*}
In \cite[p.~208]{forum-alzer}, \cite[Lemma~1.1]{batir-new}, \cite[Lemma~1.1]{batir-new-rgmia} and \cite{notes-best-simple-open, notes-best-simple.tex-rgmia}, the function
\begin{equation}\label{positivity}
[\psi'(x)]^2+\psi''(x)
\end{equation}
for $x\in(0,\infty)$ was verified to be positive. Hence, the function $[Q(x)]''>0$ and $Q(x)$ is strictly convex on $(-1,\infty)$.
\end{proof}

\begin{proof}[Proof of Theorem~\ref{Qi-Modell-thm2}]
Since the function $Q(x)$ is strictly decreasing from $(0,\infty)$ onto $\bigl(e^{-\gamma},\frac12\bigr)$, that is, the mapping $Q:(0,\infty)\mapsto\bigl(e^{-\gamma},\frac12\bigr)$ is bijective, then the proof of Theorem~\ref{Qi-Modell-thm2} is easily completed.
\end{proof}

\begin{proof}[Proof of Theorem~\ref{error-estimate}]
By the inequality~\eqref{corollary2.3-rew}, it is easy to see that
\begin{equation*}
0\le|\psi(x)-I_a(x)|\le|I_{1/2}-I_{e^{-\gamma}}|
\end{equation*}
for $x\in(c,\infty)$ and $a\in\bigl(e^{-\gamma},\frac12\bigr)$, which is equivalent to
\begin{equation*}
0\le|E_a(x)|\le\ln\frac{x+e^{-\gamma}}{x+1/2}\le\ln\frac{c+e^{-\gamma}}{c+1/2}.
\end{equation*}
The proof of Theorem~\ref{error-estimate} is complete.
\end{proof}

\begin{proof}[Proof of Theorem~\ref{harmonic-no-vu}]
Since $Q(x)$ is strictly decreasing on $(0,\infty)$, it follows that $\lim_{x\to\infty}Q(x)<Q(x)\le Q(1)$ for $x\in[1,\infty)$, which is equivalent to
\begin{gather*}
\frac12<e^{\psi(x+1)}-x\le e^{\psi(2)}-1,\\
x+\frac12<e^{\psi(x+1)}\le x+e^{1-\gamma}-1,\\
\ln\biggl(x+\frac12\biggr)<\psi(x+1)\le \ln\bigl(x+e^{1-\gamma}-1\bigr).
\end{gather*}
Taking $x=n\in\mathbb{N}$ and using the formula~\eqref{H_n-dfn} in the above inequality leads to the inequality~\eqref{new-bounds-harmonic-No}. The proof of Theorem~\ref{harmonic-no-vu} is complete.
\end{proof}

\section{Remarks}

\begin{rem}
Note that the inequality~\eqref{corollary2.3-rew} in Theorem~\ref{Qi-Modell-thm1} refines the double inequality in Corollary~2.3 of \cite{Muqattash-Yahdi}.
\end{rem}

\begin{rem}
The conclusions in Theorem~\ref{Qi-Modell-thm2} clarify the uncertain claims between lines 6\nobreakdash--7 in \cite[p.~1332]{Muqattash-Yahdi}.
\end{rem}

\begin{rem}
In \cite{qi-cui-jmaa}, the following sharp bounds for $H_n$ were established: For $n\in\mathbb{N}$,
\begin{equation}\label{ding}
\ln n+\gamma+\frac{1}{2n+1/(1-\gamma)-2}\le H(n)<\ln n+\gamma+\frac{1}{2n+1/3}.
\end{equation}
The constants $\frac{1}{1-\gamma}-2$ and $\frac13$ are the best possible.
\par
In \cite[pp.~386\nobreakdash--387]{alzer-expo-math-2006} and \cite{property-psi-ii.tex}, alternative sharp bounds for $H_n$ were presented: For $n\in\mathbb{N}$,
\begin{equation}\label{aler-harmonic-ineq}
1+\ln\bigl(\sqrt{e}\,-1\bigr)-\ln\bigl(e^{1/(n+1)}-1\bigr)\le H_n<\gamma-\ln\bigl(e^{1/(n+1)}-1\bigr).
\end{equation}
The constants $1+\ln\bigl(\sqrt{e}\,-1\bigr)$ and $\gamma$ in~\eqref{aler-harmonic-ineq} are the best possible.
\par
There have been a lot of literature devoted to bounding harmonic numbers $H_n$. For more information on $H_n$, please refer to \cite{alzer-expo-math-2006, property-psi-ii.tex, qi-cui-jmaa} and related references therein.
\par
Now we compare analytically the bounds among~\eqref{new-bounds-harmonic-No}, \eqref{ding} and~\eqref{aler-harmonic-ineq}. For this purpose, let
\begin{equation*}
f(x)=\ln\biggl(x+\frac12\biggr)+\gamma -\bigl[1+\ln\bigl(\sqrt{e}\,-1\bigr)-\ln\bigl(e^{1/(x+1)}-1\bigr)\bigr]
\end{equation*}
and
\begin{equation*}
g(x)=\ln\bigl(x+e^{1-\gamma}-1\bigr)+\ln\bigl(e^{1/(x+1)}-1\bigr).
\end{equation*}
for $x\ge1$. The function $f(x)$ may be rearranged as
\begin{equation*}
f(x)=\ln\biggl[\biggl(x+\frac12\biggr)\bigl(e^{1/(x+1)}-1\bigr)\biggr]
+\bigl[\gamma -1-\ln\bigl(\sqrt{e}\,-1\bigr)\bigr]\triangleq\ln\frac{f_1(x)}2+A,
\end{equation*}
where
\begin{gather*}
\begin{split}
f_1'(x)&=\frac{\bigl(2x^2+2x+1\bigr)\bigl[e^{1/(x+1)}-2(x+1)^2/\bigl(2x^2+2x+1\bigr)\bigr]}{(x+1)^2}\\ &\triangleq \frac{\bigl(2x^2+2x+1\bigr)f_2(1/(x+1))}{(x+1)^2},
\end{split}\\
f_2(u)=e^{u}-\frac{2}{1+(1-u)^2}=\frac{[1+(1-u)^2]e^u-2}{1+(1-u)^2} \triangleq\frac{f_3(u)}{1+(1-u)^2}
\end{gather*}
and $f_3'(u)=u^2e^u>0$ for $u\in\bigl(0,\frac12\bigr]$. Since $f_3(u)$ is increasing on $\bigl(0,\frac12\bigr]$ and $\lim_{u\to0^+}f_3(u)=0$, it follows that $f_3(u)>0$ and $f_2(u)>0$ on $\bigl(0,\frac12\bigr]$. This leads to $f_1'(x)>0$ and $f_1(x)$ being strictly increasing for $x\ge1$. Therefore, the function $f(x)$ is increasing for $x\ge1$, with
\begin{equation*}
f(1)=\ln\bigl[\bigl(\sqrt{e}\,-1\bigr)/2\bigr]+\gamma-1=-0.01731922699030\dotsm
\end{equation*}
and
\begin{align*}
\lim_{x\to\infty}f(x)&=\lim_{x\to\infty}\ln\biggl[\frac{x+1/2}{x+1} \cdot\frac{e^{1/(x+1)}-1}{1/(x+1)}\biggr]+\gamma -1-\ln\bigl(\sqrt{e}\,-1\bigr)\\
&=\gamma -1-\ln\bigl(\sqrt{e}\,-1\bigr)\\
&=0.00996779446872\dotsm.
\end{align*}
As a result, the left-hand side inequalities~\eqref{new-bounds-harmonic-No} and~\eqref{aler-harmonic-ineq} are not included each other.
\par
Numerical calculation gives
\begin{equation*}
f(2)=-0.001061745178\dotsm\quad\text{and}\quad f(3)=0.004039213518\dotsm,
\end{equation*}
thus, as $n\ge3$, the lower bound in~\eqref{new-bounds-harmonic-No} is better than the one in~\eqref{aler-harmonic-ineq}.
\par
Numerical computation gives
\begin{equation*}
g(2)=-0.00060205073286\dotsm\quad\text{and}\quad g(3)=0.00153070402207\dotsm.
\end{equation*}
It is easy to see that
\begin{equation*}
\lim_{x\to\infty}g(x)=\ln\biggl[\frac{x+e^{1-\gamma}-1}{x+1} \cdot\frac{e^{1/(x+1)}-1}{1/(x+1)}\biggr]=0.
\end{equation*}
Combining with the graph of $g(x)$ on $(1,29)$, see Figure~\ref{Inf-Fam-Dig-eps.eps},
\begin{figure}[htbp]
  \includegraphics[width=0.7\textwidth]{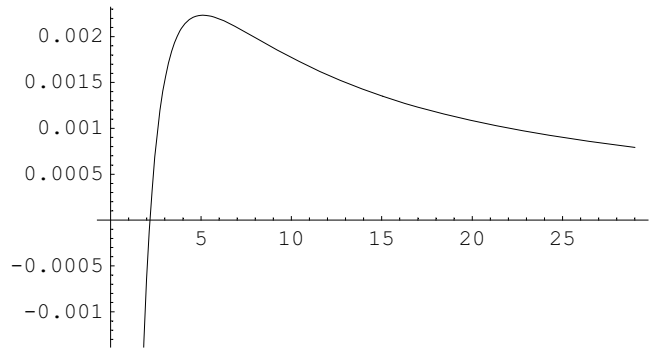}\\
  \caption{The graph of $g(x)$ on the interval $(1,29)$}\label{Inf-Fam-Dig-eps.eps}
\end{figure}
it can be possibly claimed that, as $n\ge3$, the upper bound in~\eqref{new-bounds-harmonic-No} is not better than the one in~\eqref{aler-harmonic-ineq}.
\par
Similarly, it may be showed that the upper bound in~\eqref{new-bounds-harmonic-No} is not better than the one in~\eqref{ding} as $n\ge2$, that the lower bound in~\eqref{new-bounds-harmonic-No} is not better than the one in~\eqref{ding} at all, and that the inequality~\eqref{aler-harmonic-ineq} is not better than the one~\eqref{ding} at all. Thus, it may be said that the inequality~\eqref{ding} is the best among inequalities~\eqref{new-bounds-harmonic-No}, \eqref{ding} and~\eqref{aler-harmonic-ineq}.
\end{rem}

\begin{rem}
In \cite{notes-best-simple-open, notes-best-simple.tex-rgmia}, the positivity of the function~\eqref{positivity} was generalized to a much general result including, as a particular case, its being completely monotonic on $(0,\infty)$.
\end{rem}

\begin{rem}
It is conjectured that the function $Q(x)$ defined in~\eqref{Q(x)-dfn-2} is completely monotonic on $(-1,\infty)$, that is, $(-1)^k[Q(x)]^{(k)}\ge0$ on $(-1,\infty)$ for $k\ge0$.
\end{rem}

\subsubsection*{Acknowledgements}
This manuscript was completed during the first author's visit to the \href{http://www.staff.vu.edu.au/rgmia}{RGMIA}, Victoria University, Australia between March 2008 and February 2009. The first author express his sincere appreciations on local colleagues at the \href{http://www.staff.vu.edu.au/rgmia}{RGMIA} for their invitation and hospitality throughout this period.


\begin{thebibliography}{99}

\bibitem{abram}
M. Abramowitz and I. A. Stegun (Eds), \textit{Handbook of Mathematical Functions with Formulas, Graphs, and Mathematical Tables}, National Bureau of Standards, Applied Mathematics Series \textbf{55}, 9th printing, Washington, 1970.

\bibitem{forum-alzer}
H. Alzer, \textit{Sharp inequalities for the digamma and polygamma functions}, Forum Math. \textbf{16} (2004), 181\nobreakdash--221.

\bibitem{alzer-expo-math-2006}
H. Alzer, \textit{Sharp inequalities for the harmonic numbers}, Expo. Math. \textbf{24} (2006), no.~4, 385\nobreakdash--388.

\bibitem{batir-new}
N. Batir, \textit{Some new inequalities for gamma and polygamma functions}, J. Inequal. Pure Appl. Math. \textbf{6} (2005), no.~4, Art.~103; Available online at \url{http://jipam.vu.edu.au/article.php?sid=577}.

\bibitem{batir-new-rgmia}
N. Batir, \textit{Some new inequalities for gamma and polygamma functions}, RGMIA Res. Rep. Coll. \textbf{7} (2004), no.~3, Art.~1; Available online at \url{http://www.staff.vu.edu.au/rgmia/v7n3.asp}.

\bibitem{Muqattash-Yahdi}
I. Muqattash and M. Yahdi, \textit{Infinite family of approximations of the Digamma function}, Math. Comput. Modelling \textbf{43} (2006), 1329\nobreakdash--1336.

\bibitem{notes-best-simple-open}
F. Qi, \textit{A completely monotonic function involving divided differences of psi and polygamma functions and an application}, RGMIA Res. Rep. Coll. \textbf{9} (2006), no.~4, Art.~8; Available online at \url{http://www.staff.vu.edu.au/rgmia/v9n4.asp}.

\bibitem{notes-best-simple.tex-rgmia}
F. Qi, \textit{The best bounds in Kershaw's inequality and two completely monotonic functions}, RGMIA Res. Rep. Coll. \textbf{9} (2006), no.~4, Art.~2; Available online at \url{http://www.staff.vu.edu.au/rgmia/v9n4.asp}.

\bibitem{property-psi-ii.tex}
F. Qi and B.-N. Guo, \textit{A short proof of monotonicity of a function involving the psi and exponential functions}, Available online at \url{http://arxiv.org/abs/0902.2519v1}.

\bibitem{qi-cui-jmaa}
F. Qi, R.-Q. Cui, Ch.-P. Chen, and B.-N. Guo, \textit{Some completely monotonic functions involving polygamma functions and an application}, J. Math. Anal. Appl. \textbf{310} (2005), no.~1, 303\nobreakdash--308.

\end{thebibliography}
\end{document}